\documentclass[11pt]{article}

\usepackage{latexsym}
\usepackage{amssymb}

\newtheorem{theorem}{Theorem}

\begin{document}
{
\begin{center}
{\Large\bf
On orthogonal Laurent polynomials related to the partial sums of power series.}
\end{center}
\begin{center}
{\bf S.M. Zagorodnyuk}
\end{center}

\noindent
\textbf{Abstract.}
Let $f(z) = \sum_{k=0}^\infty d_k z^k$, $d_k\in\mathbb{C}\backslash\{ 0 \}$, $d_0=1$, be a power series with a non-zero
radius of convergence $\rho$: $0 <\rho \leq +\infty$. Denote by $f_n(z)$ the n-th partial sum of $f$, and
$R_{2n}(z) = \frac{ f_{2n}(z) }{ z^n }$, $R_{2n+1}(z) = \frac{ f_{2n+1}(z) }{ z^{n+1} }$, $n=0,1,2,...$.
By the result of Hendriksen and Van Rossum there exists a linear functional $\mathbf{L}$ on Laurent polynomials, such that
$\mathbf{L}(R_n R_m) = 0$, when $n\not= m$, while $\mathbf{L}(R_n^2)\not= 0$. We present an explicit integral representation for $\mathbf{L}$
in the above case of the partial sums. We use methods from the theory of generating functions. The case of finite systems
of such Laurent polynomials is studied as well.



\section{Introduction.}

The theories of orthogonal polynomials on the real line (OPRL) and on the unit circle (OPUC) have a lot of various contributions
and applications~\cite{cit_50000_Gabor_Szego}, \cite{cit_5000_Ismail}, \cite{cit_48000_Simon_1}, \cite{cit_48000_Simon_2}.
One of their possible generalizations is a relatively new theory of biorthogonal rational functions, 
see~\cite{cit_5500_Jones_Thron_N__1984}, \cite{cit_1000_H_v_R__1986}, \cite{cit_5100_Ismail_Masson_1995},
\cite{cit_100000_Zhedanov_1999}, \cite{cit_500_Book_Bultheel___1999}, \cite{cit_400_Beckermann_Dereviagin_Zhedanov___2010} and
references therein.

Denote by $\mathcal{A}$ a set of all (formal) Laurent polynomials of the following form:
\begin{equation}
\label{l1_4}
\lambda_p x^p + \lambda_{p+1} x^{p+1} + \lambda_{p+q} x^{p+q},\qquad     p\in\mathbb{Z};\ q\in\mathbb{Z}_+,
\end{equation}
where $\lambda_j$ are complex coefficients and $x$ is an indeterminate.
Let $\{ Q_k \}_{k=0}^\infty$ be a sequence of Laurent polynomials of the following form:
\begin{equation}
\label{l1_7}
Q_{2n}(x) = \sum_{j=-n}^n \alpha_j^{(2n)} x^j,\quad \alpha_n^{(2n)}\not= 0,
\end{equation}
\begin{equation}
\label{l1_8}
Q_{2n+1}(x) = \sum_{j=-n-1}^n \alpha_j^{(2n+1)} x^j,\quad \alpha_{-n-1}^{(2n+1)}\not= 0,\qquad n\in\mathbb{Z}_+.
\end{equation}
By Proposition~1 in~\cite{cit_1000_H_v_R__1986}, if the sequence $\{ Q_k \}_{k=0}^\infty$ satisfies
\begin{equation}
\label{l1_9}
Q_{2n+1}(x) = (x^{-1} + g_{2n+1}) Q_{2n}(x) + f_{2n+1} Q_{2n-1}(x),
\end{equation}
\begin{equation}
\label{l1_10}
Q_{2n+2}(x) = (1 + g_{2n+2} x) Q_{2n+1}(x) + f_{2n+2} Q_{2n}(x),
\end{equation}
with $f_{2n+1}\not= 0$, $f_{2n+2}\not= 0$, ($n\in\mathbb{Z}_+$), and
\begin{equation}
\label{l1_12}
Q_{-1}(x) = 0,\quad Q_0(x) = \alpha_0^{(0)},
\end{equation}
then there exists a linear functional $L$: $\mathcal{A} \mapsto \mathbb{C}$, with
$L(1) = 1$, such that
\begin{equation}
\label{l1_14}
L(Q_k(x) Q_n(x))
\left\{
\begin{array}{cc} =0, & k\not=n\\
\not=0, & k=n,
\end{array}
\right.
\qquad k,n\in\mathbb{Z}_+.
\end{equation}

Recall that a $R_I$-type continuos fraction is associated with a system of monic polynomials $\{ P_n(z) \}_{n=0}^\infty$, generated by
(\cite[p. 5]{cit_5100_Ismail_Masson_1995})
\begin{equation}
\label{f1_9} 
P_n(z) = (z - \mathbf{c}_n) P_{n-1}(z) - \lambda_n (z - \mathbf{a}_n) P_{n-2}(z),\qquad n=1,2,...,
\end{equation}
where $P_{-1}(z) := 0$, $P_0(z) := 1$, and
\begin{equation}
\label{f1_10} 
\lambda_{n+1}\not= 0,\quad P_n(\mathbf{a}_{n+1})\not= 0.
\end{equation}
Polynomials $\{ P_n(z) \}_{n=0}^\infty$ are related to biorthogonal rational functions~\cite[Theorem 2.1]{cit_5100_Ismail_Masson_1995}.
The case $\mathbf{a}_n = 0$, $n\geq 2$, is related to general $T$-fractions and to the above orthogonal
Laurent polynomials $R_n$~\cite{cit_1000_H_v_R__1986}.
In fact, given a system of monic polynomials $\{ P_n(z) \}_{n=0}^\infty$, generated by~(\ref{f1_9}),(\ref{f1_10}) with $\mathbf{a}_n \equiv 0$,
and $\mathbf{c}_n \in\mathbb{C}\backslash\{ 0 \}$, $n\in\mathbb{N}$,
one can define
\begin{equation}
\label{f1_20} 
\widetilde Q_{2n}(x) = \frac{1}{ \xi_{2n} } \frac{ P_{2n}(x) }{ x^n },\ 
\widetilde Q_{2n+1}(x) = \frac{1}{ \xi_{2n+1} } \frac{ P_{2n+1}(x) }{ x^{n+1} },\ n\in\mathbb{Z}_+,
\end{equation}
where
\begin{equation}
\label{f1_22} 
\xi_k = (-1)^k \prod_{j=0}^k \mathbf{c}_j,\ \mathbf{c}_0:=1,\qquad k\in\mathbb{Z}_+.
\end{equation}
Then $\{ \widetilde Q_k \}_{k=0}^\infty$ satisfy recurrence relations~(\ref{l1_9}),(\ref{l1_10}) with
\begin{equation}
\label{f1_24} 
g_k = -\frac{1}{ \mathbf{c}_k },\ f_k = -\frac{\lambda_k \xi_{k-2}}{ \xi_k },\qquad k\in\mathbb{N}.
\end{equation}

Let 
\begin{equation}
\label{f1_29} 
f(z) = \sum_{k=0}^\infty d_k z^k,\qquad d_k\in\mathbb{C}\backslash\{ 0 \},\quad d_0=1,
\end{equation}
be a power series with a non-zero radius of convergence $\rho$: $0 <\rho \leq +\infty$. 
Denote by $f_n(z)$ the n-th partial sum of $f$, and 
\begin{equation}
\label{f1_30} 
F_n(z) = \frac{1}{d_n} f_n(z),\qquad n\in\mathbb{Z}_+.
\end{equation}
In~\cite{cit_90000_Z} it was shown that polynomials $\{ F_n(z) \}_{n=0}^\infty$ satisfy relations~(\ref{f1_9}),(\ref{f1_10})
with $\mathbf{a}_n \equiv 0$, and
\begin{equation}
\label{f1_31} 
\mathbf{c}_n = -\frac{d_{n-1}}{d_n},\quad n\in\mathbb{N};\quad  \lambda_n = \frac{ d_{n-2} }{ d_{n-1} },\quad n\geq 2.
\end{equation}
The associated orthogonal Laurent polynomials are given as follows:
\begin{equation}
\label{f1_32} 
R_{2n}(z) = \frac{ f_{2n}(z) }{ z^n },\quad R_{2n+1}(z) = \frac{ f_{2n+1}(z) }{ z^{n+1} },\qquad n\in\mathbb{Z}_+.
\end{equation}
By the above mentioned result of Hendriksen and Van Rossum there exists a linear functional $\mathbf{L}$ on Laurent polynomials, such that
$\mathbf{L}(R_n R_m) = 0$, when $n\not= m$, while $\mathbf{L}(R_n^2)\not= 0$. 
Our main purpose here is to obtain an explicit integral representation for $\mathbf{L}$.
We shall derive a generating function for $\{ R_k \}_{k=0}^\infty$, and use some methods from the theory of generating functions. 
The case of finite systems of Laurent polynomials $\{ Q_k \}_{k=0}^{2\mathbf{n}}$ ($\mathbf{n}\in\mathbb{N}$),
satisfying relations~(\ref{l1_9}),(\ref{l1_10}),
is treated in another way.
We use results on moment problems to obtain an explicit integral representation for the corresponding linear functional $L$.
The latter is done not only for the case of the partial sums, but for arbitrary such finite systems
satisfying an additional condition (see~(\ref{f2_51})).

\noindent
{\bf Notations. }
Besides the given above notations, we shall use the following definitions.
As usual, we denote by $\mathbb{R}, \mathbb{C}, \mathbb{N}, \mathbb{Z}, \mathbb{Z}_+$,
the sets of real numbers, complex numbers, positive integers, integers and non-negative integers,
respectively. 
For $k,l\in\mathbb{Z}$, we set $\mathbb{Z}_{k,l} := \{ j\in\mathbb{Z}: k\leq j\leq l\}$.
Set $\mathbb{T} := \{ z\in\mathbb{C}:\ |z|=1 \}$,
$\mathbb{D} := \{ z\in\mathbb{C}:\ |z|<1 \}$. 
By $\mathfrak{B}(\mathbb{C})$ we mean the set of all Borel subsets of $\mathbb{C}$.
By $\mathbb{P}$ we denote the set of all polynomials with complex coefficients.

\section{The partial sums of a power series and orthogonal Laurent polynomials.}

Consider a power series $f(z)$, as in~(\ref{f1_29}), its partial sums $f_n(z)$, and the associated Laurent polynomials
$\{ R_k \}_{k=0}^\infty$ from~(\ref{f1_32}).
There are various methods for deriving candidates for generating functions, as described in details in a book of McBride~\cite{cit_5700_McBride__1971},
see also a book of Rainville~\cite{cit_5150_Rainville}.
A powerful tool is a formal series manipulation, see~\cite[p. 11]{cit_5700_McBride__1971} for examples.
Interchanging the order of summation 
and using Lemma 10 on pages 56, 57 in~\cite{cit_5150_Rainville},
we can formally write:
$$ \sum_{n=0}^\infty f_n(x) t^n = \sum_{n=0}^\infty \sum_{k=0}^n d_k x^k t^n =
\sum_{n,k=0}^\infty d_k x^k t^{n+k} = $$
$$ = \sum_{n=0}^\infty t^n \sum_{k=0}^\infty d_k (xt)^k = \frac{1}{(1-t)} f(xt),\qquad |t|<1. $$
It remains to verify that the obtained candidate is valid. In the case of a generating function for $R_n$ we need some additional work, since
the correspondence between $R_n$ and $f_n$ depends on the parity of $n$.

\begin{theorem}
\label{t2_1}
Let $f(z)$ be a power series from~(\ref{f1_29}) with a non-zero radius of convergence $\rho$ ($\leq +\infty$), and $f_n(z)$ be its $n$-th partial sum.
Let $\{ R_k(z) \}_{k=0}^\infty$ be defined by~(\ref{f1_32}). 
Then the following relations hold:
\begin{equation}
\label{f2_5} 
\frac{1}{(1-t)} f(xt) = \sum_{n=0}^\infty f_n(x) t^n,\qquad t,x\in\mathbb{C}:\ |t|<1,\ |x|<\rho.
\end{equation}

$$ \left( \frac{ \sqrt{x} + 1 }{ \sqrt{x} - z } \right) f(\sqrt{x} z) + 
\left( \frac{ \sqrt{x} - 1 }{ \sqrt{x} + z } \right) f(-\sqrt{x} z) = $$
\begin{equation}
\label{f2_7} 
= 2 \sum_{n=0}^\infty R_n(x) z^n,\qquad x,z\in\mathbb{C}:\ 0<|x|<\rho,\ |z|<|\sqrt{x}|.
\end{equation}
Here by $\sqrt{x}$ we mean an arbitrarily chosen and fixed value of the square root for each $x$ (and the corresponding
values need not to form an analytic branch).

\end{theorem}
\noindent
\textbf{Proof.}
Choose an arbitrary $x$: $|x|<\rho$.
The left-hand side of~(\ref{f2_5}) is an analytic function of $t$ in $\mathbb{D}$.
Calculating its Taylor coefficient and using the Leibniz rule for the derivatives we derive relation~(\ref{f2_5}).

Let us check relation~(\ref{f2_7}).
Choose an arbitrary $x$: $0<|x|<\rho$, $t\in\mathbb{D}$, and fix an arbitrary value of $\sqrt{x}$. Denote $z = t\sqrt{x}$.
By the already established relation~(\ref{f2_5}) we may write:
$$ \left( \frac{ \sqrt{x} }{ \sqrt{x} - z } \right) f(\sqrt{x} z) = \sum_{n=0}^\infty f_n(x) (z / \sqrt{x})^n = $$
\begin{equation}
\label{f2_9} 
=  \sum_{k=0}^\infty R_{2k}(x) z^{2k} + \sqrt{x} \sum_{k=0}^\infty R_{2k+1}(x) z^{2k+1}.
\end{equation}
Denote
\begin{equation}
\label{f2_15} 
\varphi(z;x) := \left( \frac{ \sqrt{x} }{ \sqrt{x} - z } \right) f(\sqrt{x} z).
\end{equation}
Then
$$ \varphi(z;x) + \varphi(-z;x) + (\varphi(z;x) - \varphi(-z;x))/\sqrt{x} = 2\sum_{n=0}^\infty R_n(x) z^n. $$
Collecting terms with $f(\sqrt{x} z)$ and $f(-\sqrt{x} z)$ we obtain relation~(\ref{f2_7}).
$\Box$

The obtained generating functions can find various applications. For example, by the integral formula for the Taylor coefficients
one can write:
$$ R_n(x) = $$
$$ = \frac{1}{4\pi i} \oint_{|z| = |\sqrt x|/2}  
\left(
\left( \frac{ \sqrt{x} + 1 }{ \sqrt{x} - z } \right) f(\sqrt{x} z) + 
\left( \frac{ \sqrt{x} - 1 }{ \sqrt{x} + z } \right) f(-\sqrt{x} z)
\right) z^{-n-1} dz, $$
\begin{equation}
\label{f2_17} 
x \in\mathbb{C}:\ 0<|x|<\rho.
\end{equation}

Our next purpose is to obtain an explicit integral representation for the functional $\mathbf{L}$ which was discussed in
the Introduction.
In order to find a suitable candidate for a measure of integration we shall use generating functions.
It is well known that they help to establish orthogonality relations, see Section~19.3 in~\cite{cit_700_Book_Bateman___v_3}.

Observe that that the linear functional $\mathbf{L}$ is uniquely determined by the following relations:
\begin{equation}
\label{f2_19} 
\mathbf{L}(R_n) = \delta_{n,1},\qquad         n\in\mathbf{Z}_+.
\end{equation}
In fact, the span of functions $\{ R_n \}_{n=0}^\infty$ coincides with $\mathcal{A}$. If we formally apply $\mathbf{L}$ to the both
sides of relation~(\ref{f2_7}), then we obtain
\begin{equation}
\label{f2_22} 
\mathbf{L}_x(I_4(x;z)) = 2,
\end{equation}
where the superscript $x$ means that $\mathbf{L}$ acts in variable $x$, and
\begin{equation}
\label{f2_24} 
I_4(x;z) := \left( \frac{ \sqrt{x} + 1 }{ \sqrt{x} - z } \right) f(\sqrt{x} z) + 
\left( \frac{ \sqrt{x} - 1 }{ \sqrt{x} + z } \right) f(-\sqrt{x} z).
\end{equation}
It is convenient to introduce a new variable $y=\sqrt{x}$, and write
\begin{equation}
\label{f2_27} 
I_4(y^2;z) := \left( \frac{ y + 1 }{ y - z } \right) f(y z) + 
\left( \frac{ y - 1 }{ y + z } \right) f(-y z).
\end{equation}

We can multiply the right-hand side of~(\ref{f2_27}) by an arbitrarily chosen function $a(y)$, and then calculate
some contour integrals $\oint I_4(y^2;z) a(y) dy$, trying to obtain the value $2$. 
In this manner we obtain a candidate which is described in the next theorem.

\begin{theorem}
\label{t2_2}
Let $f(z)$ be a power series as in~(\ref{f1_29}) with a non-zero radius of convergence $\rho$, and $f_n(z)$ be its $n$-th partial sum.
Define $\{ R_k(z) \}_{k=0}^\infty$ by relations~(\ref{f1_32}). Let $\mathbf{L}$ be a linear functional on $\mathcal{A}$,
such that $\mathbf{L}(R_n R_m) = 0$, when $n\not= m$, while $\mathbf{L}(R_n^2)\not= 0$, $n,m\in\mathbb{Z}_+$. 
Then the following integral representation holds:
\begin{equation}
\label{f2_29} 
\mathbf{L}(R) = \frac{1}{2\pi i} \oint_{ |y|=c } R(y^2) \frac{dy}{ yf(y^2) },\qquad        R\in\mathcal{A}.    
\end{equation}
Here $c$ is an arbitrary positive number which is less than $\sqrt{\rho}$, and less than $\widehat\rho$, where $\widehat\rho$ is the radius of convergence
for the Maclaurin series of $1/f(y^2)$.~\footnote{Since $f(0)=1$, it follows by continuity that $\widehat\rho > 0$.}

\end{theorem}
\noindent
\textbf{Proof.}
Denote the right-hand side of~(\ref{f2_29}) by $\mathcal{L}(R)$.
Let us check that $\mathcal{L}$ has property~(\ref{f2_19}).
Choose an arbitrary $n\in\mathbb{Z}_+$. 
By~(\ref{f1_32}) we may write:
$$ \mathcal{L}(R_{2n}) = \frac{1}{2\pi i} \oint_{ |y|=c } f_{2n}(y^2) \frac{dy}{ y^{2n+1} f(y^2) } = $$
$$ = \frac{1}{2\pi i} \oint_{ |y|=c } \left( f(y^2) - \sum_{j=2n+1}^\infty d_j y^{2j} \right)
\frac{dy}{ y^{2n+1} f(y^2) } = $$
$$ = \frac{1}{2\pi i} \oint_{ |y|=c } \left( 1 - \left(\sum_{j=2n+1}^\infty d_j y^{2j} \right) ( f(y^2) )^{-1} \right)
\frac{dy}{ y^{2n+1} } = $$
$$ = \frac{1}{ (2n)! } \left[ 1 - \left(\sum_{j=2n+1}^\infty d_j y^{2j} \right) ( f(y^2) )^{-1} \right]^{(2n)} (0) = 
\left\{ \begin{array}{cc} 1, & n=0\\
              0, &  n>0            \end{array}
\right. . $$
The last equality can be justified, for example, using the Leibniz rule for derivatives.
In a similar way, we may write
$$ \mathcal{L}(R_{2n+1}) = \frac{1}{2\pi i} \oint_{ |y|=c } f_{2n+1}(y^2) \frac{dy}{ y^{2n+3} f(y^2) } = $$
$$ = \frac{1}{2\pi i} \oint_{ |y|=c } \left( f(y^2) - \sum_{j=2n+2}^\infty d_j y^{2j} \right)
\frac{dy}{ y^{2n+3} f(y^2) } = $$
$$ = \frac{1}{ (2n+2)! } \left[ 1 - \left(\sum_{j=2n+2}^\infty d_j y^{2j} \right) ( f(y^2) )^{-1} \right]^{(2n+2)} (0) = 0. $$ 
This completes the proof. $\Box$

For example, consider the following function:
\begin{equation}
\label{f2_35} 
f(z) = e^{bz} \prod_{j=0}^m (1 - a_j z)^{-\lambda_j},\qquad b\geq 0,\ 0<a_j<1,\ \lambda_j>0;\ m\in\mathbb{N}.
\end{equation}
The corresponding Maclaurin series converges in $K := \{ z\in\mathbb{C}:\ |z| < \min(1/a_0,...,1/a_m) \}$, which contains $\mathbb{D}$.
Observe that this series has positive coefficients and the constant term $f(0)=1$.
Denote $\rho := \min(1/a_0,...,1/a_m)$.
The function
$$ 1/ f(y^2) = e^{-b y^2} \prod_{j=0}^m (1 - a_j y^2)^{\lambda_j}, $$
is analytic in $K_1 := \{ z\in\mathbb{C}:\ |z| < \sqrt{ \rho } \}$. Therefore, in this case we may 
apply Theorem~\ref{t2_2} and write
\begin{equation}
\label{f2_39} 
\mathbf{L}(R) = \frac{1}{2\pi i} \oint_{ \mathbb{T} } R(y^2) e^{-b y^2} \prod_{j=0}^m (1 - a_j y^2)^{\lambda_j}
\frac{dy}{ y },\qquad        R\in\mathcal{A}.    
\end{equation}

Notice that one may also use various cases of the generalized hypergeometric function as $f(z)$, and investigate the corresponding
partial sums and orthogonal Laurent polynomials.

\noindent
\textbf{Finite systems of Laurent polynomials.}
Fix an arbitrary $\textbf{n}\in\mathbb{N}$. Denote by $\mathcal{A}_{\textbf{n}}$
a set of all (formal) Laurent polynomials of the following form:
\begin{equation}
\label{f2_48}
\sum_{j=-\textbf{n}}^{\textbf{n}} \lambda_j x^j,\qquad     \lambda_j\in\mathbb{C},
\end{equation}
where $x$ is an indeterminate.
Let $\{ Q_k \}_{k=0}^{2\mathbf{n}}$ be a set of Laurent polynomials which have forms as in~(\ref{l1_7}),(\ref{l1_8}).
Of course, these polynomials belong to $\mathcal{A}_{\textbf{n}}$ and span it.

Suppose that $\{ Q_k \}_{k=0}^{2\mathbf{n}}$ satisfy relations~(\ref{l1_9}),(\ref{l1_10}) with
some $f_{2n+1}\not= 0$, $f_{2n+2}\not= 0$, ($n\in\mathbb{Z}_{0,\mathbf{n}-1}$), and $Q_{-1}(x) := 0$.
We can expand the sequences of complex coefficients $\{ f_k \}_{k=0}^{2\mathbf{n}}$, $\{ g_k \}_{k=0}^{2\mathbf{n}}$,
to infinite complex sequences $\{ f_k \}_{k=0}^\infty$, $\{ g_k \}_{k=0}^\infty$, $f_k\not=0$, in an arbitrary way.
Then we extend the sequence of $Q_k$ by relations~(\ref{l1_9}),(\ref{l1_10}) to a sequence $\{ Q_k \}_{k=0}^\infty$. 
By the result of Hendriksen and Van Rossum there exists a linear functional $L$ on Laurent polynomials, such that
$L(Q_n Q_m) = 0$, when $n\not= m$, while $Q(L_n^2)\not= 0$.

Denote by $L_{\textbf{n}}$ the restriction of $L$ to the set $\mathcal{A}_{\textbf{n}}$. We have
\begin{equation}
\label{f2_50}
L_{\textbf{n}} (Q_k Q_n) 
\left\{
\begin{array}{cc} =0, & k\not=n\\
\not=0, & k=n,
\end{array}
\right.
\qquad k,n\in\mathbb{Z}_{0,2\mathbf{n}}.
\end{equation}

\begin{theorem}
\label{t2_3}
Let $\{ Q_k \}_{k=0}^{2\mathbf{n}}$ be a finite set of Laurent polynomials, having forms as in~(\ref{l1_7}),(\ref{l1_8});
$\textbf{n}\in\mathbb{N}$.
Suppose that $\{ Q_k \}_{k=0}^{2\mathbf{n}}$ satisfy relations~(\ref{l1_9}),(\ref{l1_10}) with
some $f_{2n+1}\not= 0$, $f_{2n+2}\not= 0$, ($n\in\mathbb{Z}_{0,\mathbf{n}-1}$), and $Q_{-1}(x) := 0$.
Suppose that the corresponding linear functional $L_{\textbf{n}}$, having property~(\ref{f2_50}),
satisfies the following condition:
\begin{equation}
\label{f2_51}
a := L_{\textbf{n}} (x^{-\mathbf{n}}) \not= 0,
\end{equation}
Then $L_{\textbf{n}}$ admits the following integral representation:
\begin{equation}
\label{f2_52}
L_{\textbf{n}} (Q) = 
\int Q(z) a z^\mathbf{n} d\mu,\qquad Q\in\mathcal{A}_{\textbf{n}},
\end{equation}
where $\mu$ is a finitely atomic positive measure on $\mathfrak{B}(\mathbb{C})$.

\end{theorem}
\noindent
\textbf{Proof.}
Denote
$$ s_k = \frac{1}{a} L_{\textbf{n}} (x^{k-\mathbf{n}}),\qquad    k=0,1,...,2\mathbf{n}. $$
Consider the following moment problem (see~\cite{cit_3700_Z}, \cite{cit_95000_Zagorodnyuk}):
find a (non-negative) measure $\mu$ on $\mathfrak{B}(\mathbb{C})$
such that
\begin{equation}
\label{f2_54}
\int_{\mathbb{C}} z^k d\mu(z) = s_k,\qquad  k\in\mathbb{Z}_{0,2\mathbf{n}}.
\end{equation}
Since $s_0=1$, it is solvable and according to Algorithm~1 in~\cite{cit_95000_Zagorodnyuk}
it has a finitely atomic solution.
For an arbitrary $Q\in\mathcal{A}_{\mathbf{n}}$ we may write:
$$  Q(z) = \sum_{k=0}^{2\mathbf{n}} a_k z^{k-\mathbf{n}},\qquad a_k\in\mathbb{C}. $$
Substituting this expression into both sides of relation~(\ref{f2_52}), we shall obtain the same value.
The proof is complete.
$\Box$

We remark that Algorithm~1 in~\cite{cit_95000_Zagorodnyuk} provides an explicit procedure for the construction
of the corresponding atomic solution (see Example~1 in~\cite{cit_95000_Zagorodnyuk}).

Thus, for large classes of orthogonal Laurent polynomials it was here possible to construct explicit measures and
integral representations. The theories of generating functions and moment problems have shown their power and usefulness as
the corresponding tools of the investigation.

\vspace{1.5cm} 

V. N. Karazin Kharkiv National University \newline\indent
School of Mathematics and Computer Sciences \newline\indent
Department of Higher Mathematics and Informatics \newline\indent
Svobody Square 4, 61022, Kharkiv, Ukraine

Sergey.M.Zagorodnyuk@gmail.com; Sergey.M.Zagorodnyuk@univer.kharkov.ua

}
\end{document}